\documentclass{article}
\usepackage[utf8]{inputenc}
\usepackage[usenames]{color}
\usepackage{amsmath}
\usepackage{epsfig,latexsym}
\usepackage{epsf}
\usepackage{amsmath,amsfonts,amstext,amssymb,amsbsy,amsopn,amsthm,eucal}
\usepackage{latexsym}
\usepackage{epsfig,latexsym}
\usepackage[usenames]{color}
\usepackage{dirtytalk}
\newcommand{\A}{{\rm A}}
\title{Ouroboros Spaces:\\ An Intuitive Approach to Self-Referential Functional Analysis with Applications to Probability Theory}
\author{Nathan Thomas Provost$\footnote{Student of Applied Mathematics and Statistics at Brown University. University Email: nathan\_provost@brown.edu}$}
\date{}
\begin{document}
\maketitle

\begin{center}
\textbf{Abstract}
\end{center}

\small{In this paper, we aim to introduce the concept of the Ouroboros space and the complimentary concept of the Ouroboros function by using the Ouroboros equation \cite{1} as our starting point. We start with a few univariate definitions, and then move on to multivariate definitions. We contextualize and motivate these concepts with a few examples. Then, we discuss a few aspects from probability theory that are relevant to the idea of an Ouroboros space, eventually proving two critical theorems. We briefly introduce the case of mixed domains, and summarize our findings, while emphasizing the importance of self-referential functions.}

\section*{Introduction}

\normalsize Ancient history and historical artifacts have had an indescribable impact on the development of modern society. Recently, the symbol of a self-consuming snake (historically known as the \textit{Ouroboros}) has been the inspiration for various ideas in mathematics, computer science, and biology. Chiefly, it has been used to name an equation that embodies the idea of a self-reference, the likes of which was the central topic of a fascinating paper, written within the last decade, focusing on the \say{the various manifestations or ways in which [the] Ouroboros equation has emerged.} \cite{1} We instead aim to give a new, intuitive definition to the function spaces that consist of solutions to this enigmatic equation. Naturally, we begin with the univariate case for which we prove some interesting properties, but we quickly advance to the multivariate case, which allows us to tie probability theory to Ouroboros spaces in an interesting way.

\section*{Ouroboros Spaces}
We are chiefly concerned with the idea of the \textit{Ouroboros Equation} \cite{1}:
\[f(f)=f\]
We will first build our definitions using the univariate case, given similarly by $f(f(x))=f(x), \ \forall x$. Suppose that $A$ is a set and $B$ is another set. Then we will call $\textbf{\textit{O}}(A)$ the \textit{Ouroboros Function Space}, or more laconically, the \textit{Ouroboros Space} for the domain given by $A$. To be more explicit, we can deduce that:
\[f\in\textbf{\textit{O}}(A) \longrightarrow f:A\rightarrow B \ni f(f(x))=f(x), \ \ \forall x\in A\]
We call such a function an \textit{Ouroboros function}. It obviously follows that $\textbf{\textit{O}}(A)$ is the set of all functions $f$ that satisfy this condition. From this realization, we can rigorously define the Ouroboros Space like so, though we will later revise and perfect this definition:
\[\textbf{\textit{O}}(A)=\left\{f:A\rightarrow B \mid f(f(x))=f(x), \ \forall x\in A\right\}\]
To motivate this concept, we provide a brief example.\\

\noindent \textbf{Example 1}: Suppose that $f:\mathbb{R}\rightarrow\mathbb{R} \ni f(x)=x$. Let $A=\mathbb{R}=B$. Then we can now obviously write $f:A\rightarrow B \ni f(x)=x$. Let $x_0\in \mathbb{R}=A$ be an arbitrary value in the domain for $x$. By definition, $f(x_0)=x_0$, which means that $f(f(x_0))=f(x_0)=x_0, \ \forall x_0\in A=\mathbb{R}$. Therefore, it holds that $f\in\textbf{\textit{O}}(\mathbb{R})$ from our definition. A nearly identical proof holds for a complex domain. $\qedsymbol$ \\

The previous case gives an immediate and quite obvious example of an Ouroboros function, but it does so under a particular condition. Namely, the last example unfolds under the condition that $A=B$, but this seems a bit bold of an assumption to make. In the following example, we prove that $f:A\rightarrow B$ can be an Ouroboros function without the condition that $A=B$.\\

\noindent \textbf{Example 2}: Suppose we have some real number $c\in\mathbb{R}$ which determines the function $g:\mathbb{R}\rightarrow \{c\} \ni g(x)=c$. Obviously, we can let $A=\mathbb{R}$ and $B=\{c\}$, observing that $B=\{c\}\subset\mathbb{R}=A$. From this, we have a function of the form $g:A\rightarrow B$ again. If we let $x_0\in\mathbb{R}=A$ once more, we see that $g(x_0)=c$, which means that $g(g(x_0))=g(c)$, but from our definition of this constant function, $g(c)=c$, so $g(g(x_0))=g(c)=c=g(x_0), \ \forall x_0\in\mathbb{R}$. From this deduction, we ultimately conclude that $g\in\textbf{\textit{O}}(\mathbb{R})$. $\qedsymbol$ \\

We have now shown that it is possible for an Ouroboros function to exist in the form $f:A\rightarrow B$ where $B\subseteq A$. Furthermore, we note that the case where $A\nsubseteq B$ and $B\nsubseteq A$ can never hold, since $f(f(x))$ would always be undefined. A natural question after all of this is: \textit{can we have an Ouroboros function where} $A\subset B$? In short, the answer is no, which leads to the following lemma:\\

\noindent \textbf{Lemma 1}: Suppose we have a univariate function $f$. If $f\in\textbf{\textit{O}}(A)$, where $\textbf{\textit{O}}(A)=\left\{f:A\rightarrow B : f(f(x))=f(x), \ \forall x\in A\right\}$ and $f:A\rightarrow B$, then $B\subseteq A$.\\

\noindent\textbf{Proof}: We will prove this lemma by employing the method of proof by contradiction. We \textbf{assume} we have an Ouroboros function $f\in\textbf{\textit{O}}(A)$, where $f:A\rightarrow B$ and $A\subset B$, noting that $f(x)$ is defined for \textit{inputs} \textbf{only} in the domain $A$. Let us define $C=B/A$ such that $C\subset B$, but $C\nsubseteq A$ and $A\nsubseteq C$. This further implies that $A\cup C=B$. From these definitions, we deduce that $\exists \ x_0 \in A \ni f(x_0)\in C$, which since $C\subset B$, means that $f(x_0)\in B$ also. However, we stated that $C\nsubseteq A$ and $A\nsubseteq C$, which means that since $f(x_0)\in C$, it also holds that $f(x_0)\notin A$. This means that $f(f(x_0))$ is undefined because $f(x_0)\notin A$, implying that $f(f(x_0))\neq f(x_0)$. This would mean $f$ is not an Ouroboros function, which is a clear contradiction of our assumptions. Therefore, if $f\in\textbf{\textit{O}}(A)$ where $f:A\rightarrow B$, then $B\subseteq A$. $\qedsymbol$\\

This lemma gives us valuable information about Ouroboros functions and Ouroboros spaces. Naturally, we will now revise our initial definition to provide the true definition of an Ouroboros space for the domain given by $A$.\\

\noindent\textbf{Definition 1}: The \textbf{Ouroboros Space for the Domain Given by} $A$ is a function space given by:
\[\textbf{\textit{O}}(A)=\left\{f:A\rightarrow B \mid f(f(x))=f(x), \ \forall x\in A, \ \forall B\subseteq A \right\}\]
If $f\in\textbf{\textit{O}}(A)$, then $f(x)$ is said to be an \textbf{Ouroboros Function for the Domain Given by} $A$.\\

Now that we have obtained the proper definition of an Ouroboros space, we can apply example 1 to several iterations of the parent linear function. If $f:\mathbb{Z}\rightarrow\mathbb{Z}\ni f(x)=x, \ \forall x\in\mathbb{Z}$, then $f\in\textbf{\textit{O}}(\mathbb{Z})$. We could keep changing the domain and codomain of our function to find examples of functions from $\textbf{\textit{O}}(\mathbb{N})$,$\textbf{\textit{O}}(\mathbb{Q})$, and $\textbf{\textit{O}}(\mathbb{C})$, but these examples would just be boring iterations of the same idea. The next progression regarding Ouroboros spaces and Ouroboros functions comes from the concept of a multivariate extension. From this extension, as we will see, comes a relatively shocking result about a commonly used numerical value that is absolutely essential to the formal and natural sciences. To start this process, we begin with a new, rough definition of a multivariate Ouroboros space that we will later perfect as we did with our univariate definition:
\[\textbf{\textit{O}}(A^n)=\left\{f:A^n\rightarrow B \mid f(f(\textbf{x}),...,f(\textbf{x}))=f(\textbf{x}), \ \forall \textbf{x}\in A^n\right\}\]
wherein we assume $\textbf{x}=[x_1 \ ... \ x_n]^T\in A^n$ and $f(\textbf{x})=f(x_1,...,x_n)$, so that we also mean $f(f(\textbf{x}),...,f(\textbf{x}))=f(f(x_1,...,x_n),...,f(x_1,...,x_n))$. We also assume that for $i=1,...,n\in\mathbb{N}$, $x_i\in A$. We will now state the multivariate analog of lemma 1 and give a proof of its implications.\\

\textbf{Lemma 2}: Suppose we have a multivariate function $f$. If $f\in\textbf{\textit{O}}(A^n)$, where we have $\textbf{\textit{O}}(A^n)=\left\{f:A^n\rightarrow B \mid f(f(\textbf{x}),...,f(\textbf{x}))=f(\textbf{x}), \ \forall \textbf{x}\in A^n\right\}$ and $f:A^n\rightarrow B$, then $B\subseteq A$.\\

\noindent\textbf{Proof}: We again prove this lemma through proof by contradiction. Assume $f\in\textbf{\textit{O}}(A^n)$, $f:A^n\rightarrow B$, and $A\subset B$, again noting that $f$ is defined only for inputs from its domain $A$. Then we again define $C=B/A$ such that $C\subset B$, but $C\nsubseteq A$ and $A\nsubseteq C$. This further implies that $A\cup C=B$. From these definitions, we deduce that $\exists \ c_1,...,c_n \in A \ni f(c_1,...,c_n)\in C$, which since $C\subset B$, means that $f(c_1,...,c_n)\in B$ also. Since $C\nsubseteq A$ and $A\nsubseteq C$, $f(c_1,...,c_n)\in C$ means that $f(c_1,...,c_n)\notin A$. This means that $f(f(c_1,...,c_n),...,f(c_1,...,c_n))$ is undefined, since $f(c_1,...,c_n)\notin A$, implying that this value is not in the domain, which in turn means that $f(f(c_1,...,c_n),...,f(c_1,...,c_n))\neq f(c_1,...,c_n)$. This means that $f$ is not an Ouroboros function, which contradicts our initial assumption. Therefore, if $f\in\textbf{\textit{O}}(A^n)$ and $f:A^n\rightarrow B$, then $B\subseteq A$. $\qedsymbol$ \\

From this lemma, we can make our second definition. This general statement covers Ouroboros spaces and functions of an arbitrary number of variables.\\

\textbf{Definition 2}: The \textbf{\textit{General} Ouroboros Space for the Domain Given by} $A^n$ is a function space given by:
\[\textbf{\textit{O}}(A^n)=\left\{f:A^n\rightarrow B \mid f(f(\textbf{x}),...,f(\textbf{x}))=f(\textbf{x}), \ \forall \textbf{x}\in A^n, \ \forall B\subseteq A\right\}\]
If $f\in\textbf{\textit{O}}(A^n)$, then $f(\textbf{x})$ is said to be an \textbf{Ouroboros Function for the Domain Given by} $A^n$, where $\textbf{x}=[x_1 \ ... \ x_n]^T\in A^n$\\

Naturally, we might want to motivate this abstract higher level definition with a few examples. First, we consider two examples with particular dimension sizes (namely $n=2$ and $n=3$). Then, we will prove a general result that ties the ideas of an Ouroboros space and an Ouroboros function to a fundamental concept in probability theory.\\

\textbf{Example 3}: Consider the function $f:\mathbb{R}^2\rightarrow\mathbb{R}$ where $f(x,y)=\frac{1}{2}(x+y)$. We can choose any $x_0\in\mathbb{R}$ and $y_0\in\mathbb{R}$, noting that $A=\mathbb{R}=B$ which means $f:A^2\rightarrow B$. By definition, $f(x_0,y_0)=\frac{1}{2}(x_0+y_0)$, so it follows that $f(f(x_0,y_0),f(x_0,y_0))=\frac{1}{2}(\frac{1}{2}(x_0+y_0)+\frac{1}{2}(x_0+y_0))=\frac{1}{2}(x_0+y_0)=f(x_0,y_0)$. Thus, by definition, $f\in\textbf{\textit{O}}(\mathbb{R}^2)$. $\qedsymbol$\\

\textbf{Example 4}: Now, consider the function $g:\mathbb{R}^3\rightarrow\mathbb{R}$ where $g(x,y,z)=\frac{1}{3}(x+y+z)$. Again, we choose any $x_0\in\mathbb{R}$, $y_0\in\mathbb{R}$, and $z_0\in\mathbb{R}$, noting that $A=\mathbb{R}=B$ such that $g:A^3\rightarrow B$. By definition, $g(x_0,y_0,z_0)=\frac{1}{3}(x_0+y_0+z_0)$, so it naturally follows that we can go on to write $g(g(x_0,y_0,z_0),g(x_0,y_0,z_0),g(x_0,y_0,z_0))=\frac{1}{3}(\frac{1}{3}(x_0+y_0+z_0)+\frac{1}{3}(x_0+y_0+z_0)+\frac{1}{3}(x_0+y_0+z_0))=\frac{1}{3}(x_0+y_0+z_0)=g(x_0,y_0,z_0)$. Therefore, following again from our definition, we can say that $g\in\textbf{\textit{O}}(\mathbb{R}^3)$. $\qedsymbol$

\section*{Probability Theory with Ouroboros Functions}

We are starting to notice a pattern in our examples from the previous section regarding higher dimensional Ouroboros functions. In example 3, our function is essentially just the average of two numbers ($x$ and $y$), and in example 4, our function is just the average of three numbers ($x$, $y$, and $z$). This leads us to a general theorem regarding arithmetic averages and Ouroboros spaces.\\

\textbf{Theorem 1}: Suppose we have a \textbf{general arithmetic average function} of $n$ arbitrary variables $x_i\in\mathbb{R}, \ \forall i=1,...,n\in\mathbb{N}$, denoted by:
\[\A:\mathbb{R}^n\rightarrow\mathbb{R} \ \ni \A(x_1,...,x_n)=\frac{1}{n}\sum_{i=1}^{n}x_i.\]
For any arbitrary dimension $n\in\mathbb{N}$, $\A(x_1,...,x_n)\in\textbf{\textit{O}}(\mathbb{R}^n), \ \forall x_i\in\mathbb{R}$.\\

\noindent\textbf{Proof}: We first must state the general definition of $\textbf{\textit{O}}(\mathbb{R}^n)$ based on our general Ouroboros space definition for $\textbf{x}=[x_1 \ ... \ x_n]^T\in \mathbb{R}^n$:
\[\textbf{\textit{O}}(\mathbb{R}^n)=\left\{f:\mathbb{R}^n\rightarrow B \mid f(f(\textbf{x}),...,f(\textbf{x}))=f(\textbf{x}), \ \forall \textbf{x}\in \mathbb{R}^n, \ \forall B\subseteq \mathbb{R}\right\}\]
In the case of $\A(x_1,...,x_2)$, it is obvious that $B=\mathbb{R}\subseteq\mathbb{R}$, so this condition is met. Now we simply need to show that $\A$ is a self-referential function. We begin as usual by choosing an arbitrary string of $n$ values from the domain, given by $c_1,...,c_n\in\mathbb{R}$. By the definition of our function, we know that:\\
\[\A(c_1,...,c_n)=\frac{1}{n}\sum_{i=1}^{n}c_i=\frac{1}{n}(c_1+...+c_n)\]
We can now observe that the value of $\A(\A(c_1,...,c_n),...,\A(c_1,...,c_n))$ is given by:
\[\A(\A(c_1,...,c_n),...,\A(c_1,...,c_n))=\A(\frac{1}{n}(c_1+...+c_n),...,\frac{1}{n}(c_1+...+c_n))=\]
\[\frac{1}{n}\sum_{i=1}^{n}x_i, \ \forall i=1,...,n \ni x_i=\frac{1}{n}(c_1+...+c_n)\]
Knowing this, we can deduce that:
\[\A(\A(c_1,...,c_n),...,\A(c_1,...,c_n))=\frac{1}{n}(\frac{1}{n}(c_1+...+c_n)+...+\frac{1}{n}(c_1+...+c_n))=\]
\[\frac{1}{n}(\frac{n}{n}(c_1+...+c_n))=\frac{1}{n}(c_1+...+c_n)=\frac{1}{n}\sum_{i=1}^{n}c_i=\A(c_1,...,c_n)\]
Therefore, since $\A:\mathbb{R}^n\rightarrow\mathbb{R}$, $\mathbb{R}\subseteq\mathbb{R}$, and $\A(\A(x_1,...,x_n),...,\A(x_1,...,x_n))=\A(x_1,...,x_n)$, we conclude that $\A(x_1,...,x_n)\in\textbf{\textit{O}}(\mathbb{R}^n)$. $\qedsymbol$ \\

This is quite an idea, for it is difficult to think of functions that satisfy the Ouroboros equation by simply brainstorming strange functions. This relationship seems vaguely familiar to a relationship in probability theory regarding expected value. For a random variable $X$, it always holds that $\mathbb{E}[\mathbb{E}[X]]=\mathbb{E}[X]$, and this property seems to echo the Ouroboros equation. Moreover, whether we are discussing $\A(\A(\textbf{x}),...,\A(\textbf{x}))$ or $\mathbb{E}[\mathbb{E}[X]]$, we are referring to the average of averages, both of which yield a self-same average.

We now recall the Strong Law of Large Numbers, which is adapted from the definition given by Evans \cite{2}. Suppose we have independent, identically distributed, integrable random random variables $X_1,...,X_n$ on the probability space $\mathbb{S}=(\Omega, \mathbb{F}, \mathbb{P})$, where $\Omega$ is the sample space, $\mathbb{F}$ is an appropriate $\sigma$-algebra (such that each $X_i$ is $\mathbb{F}$-measurable), and $\mathbb{P}$ is a probability measure. Also, let $i=1,...,n\in\mathbb{N}$. Then:
\[\mathbb{P}\left(\lim_{n\rightarrow\infty}\frac{1}{n}\sum_{i=1}^{n}X_i=\mathbb{E}[X_i]\right)=1\]
Again, we have adapted Evans' definition with slight notational differences \cite{2}, which can be rephrased as:
\[\lim_{n\rightarrow\infty}\frac{1}{n}\sum_{i=1}^{n}X_i=\mathbb{E}[X_i] \ \rm{(almost \ surely)}\]
Assume that for each $X_i$ where $i=1,...,n\in\mathbb{N}$, we can furthermore say that $X_i:\Omega\rightarrow\mathbb{R}, \ \forall i=1,...,n$. If we use our previous definition of:
\[\A(x_1,...,x_n)=\frac{1}{n}\sum_{i=1}^{n}x_i\]
Then we can go on to write:
\[\A(X_1,...,X_n)=\frac{1}{n}\sum_{i=1}^{n}X_i\]
Furthermore, if we impose an infinite limit on this expression, we find that:
\[\lim_{n\rightarrow\infty}\A(X_1,...,X_n)=\lim_{n\rightarrow\infty}\frac{1}{n}\sum_{i=1}^{n}X_i=\mathbb{E}[X_i]\ \rm{(almost \ surely)}\]
From our general definition, we can equivalently define:
\[\textbf{\textit{O}}(\mathbb{R}^\infty)=\left\{f:\mathbb{R}^\infty\rightarrow B \mid f(f(\textbf{x}),f(\textbf{x}),...)=f(\textbf{x}), \ \forall \textbf{x}\in \mathbb{R}^\infty, \ \forall B\subseteq \mathbb{R}\right\}\]
Here, we can say $\textbf{x}$ takes the form:
\[\textbf{x}=\lim_{n\rightarrow\infty}[x_1 \ ... \ x_n]^T\]
Finally, from Theorem 1, we can make the ultimate conclusion that:
\[\lim_{n\rightarrow\infty}\A(X_1,...,X_n)=\mathbb{E}[X_i]\in\textbf{\textit{O}}(\mathbb{R}^\infty) \ \rm{(almost \ surely)}\]
We can generalize our previous discussion into the form of a second theorem that encompasses all of our assumptions and conclusions, summarized with a quick proof.\\

\textbf{Theorem 2}: Suppose $X_1,...,X_n$ is a sequence of independent, identically distributed, integrable, $\mathbb{F}$-measurable random variables defined on the probability space $\mathbb{S}=(\Omega, \mathbb{F}, \mathbb{P})$, such that $\mathbb{F}$ is an appropriate $\sigma$-algebra for $\Omega$. Moreover, assume that each $X_i(\omega)=X_i$ takes the form $X_i:\Omega\rightarrow\mathbb{R}, \ \forall \omega\in\Omega, \ \forall i=1,...,n\in\mathbb{N}$. Then, in accordance with the Strong Law of Large Numbers, for any $i=1,...,n$, it holds that:
\[\mathbb{E}[X_i]\in\textbf{\textit{O}}(\mathbb{R}^\infty) \ \rm{(almost \ surely)}\]

\noindent\textbf{Proof}: Suppose that all of the conditions in Theorem 2 are met as needed. Let the general arithmetic average function ($\A(x_1,...,x_n)$) be defined as we have previously defined it. Observe that under the previously stated conditions, $\A(X_1,...,X_n):\mathbb{R}^n\rightarrow\mathbb{R}$ for any positive integer $n$. Naturally, as we said before:
\[\lim_{n\rightarrow\infty}\A(X_1,...,X_n)=\lim_{n\rightarrow\infty}\frac{1}{n}\sum_{i=1}^{n}X_i\]
Similarly, we see that when $n\rightarrow\infty$, $\A(X_1,...,X_n):\mathbb{R}^\infty\rightarrow\mathbb{R}$, which in turn means that by our definition of the general Ouroboros space and Theorem 1:
\[\lim_{n\rightarrow\infty}\A(X_1,...,X_n)\in\textbf{\textit{O}}(\mathbb{R}^\infty)\]
By the Strong Law of Large Numbers, we can conclude that for $i=1,...,n\in\mathbb{N}$:
\[\mathbb{E}[X_i]=\lim_{n\rightarrow\infty}\frac{1}{n}\sum_{i=1}^{n}X_i=\lim_{n\rightarrow\infty}\A(X_1,...,X_n)\in\textbf{\textit{O}}(\mathbb{R}^\infty) \ \rm{(almost \ surely)} \ \qedsymbol\]

Breaking away from this fascinating result, we should now briefly discuss the case of mixed domains, or rather, domains that can be expressed as Cartesian products of unequal sets. Suppose we have a domain given by $\Delta=[A_1\times...\times A_n]$ where $A_k\neq A_j, \ \forall k\neq j$. We denote the Ouroboros space for all $i=1,...,n\in\mathbb{N}$ such that $\Delta=[A_1\times...\times A_n]$ for this domain by the following definition:
\[\textbf{\textit{O}}(\Delta)=\left\{f:\Delta\rightarrow B \mid f(f(\textbf{x}),...,f(\textbf{x}))=f(\textbf{x}), \ \forall \textbf{x}\in \Delta, \ \forall B\subseteq A_i\right\}\]
This slightly altered definition requires that $B$ is a subset of \textbf{each} $A_i$, and we will now prove the following lemma stating that $A_i\not\subset B$, for all $i$.\\

\newpage
\textbf{Lemma 3}: If $\Delta=[A_1\times...\times A_n]$, $f\in\textbf{\textit{O}}(\Delta)$, and $f:\Delta\rightarrow B$, then $B\subseteq A_i$ for all $A_i$ and $i=1,...,n$.\\

\noindent\textbf{Proof}: Assume $f:\Delta\rightarrow B$, $f\in\textbf{\textit{O}}(\Delta)$ with $\Delta=[A_1\times...\times A_n]$, and $\exists A_i\subset B$ for some $i\in{1,...,n}\subset\mathbb{N}$, in which there is some $x_i\in A_i$ such that $f(x_i)\in B/A_i=C$, which means that $f(x_1,...,x_i,...,x_n)\notin A_i$ because $A_i\nsubseteq C$ and $C \nsubseteq A_i$. Therefore, $f(f(x_1,...,x_i,...,x_n),...,f(x_1,...,x_i,...,x_n))\neq f(x_1,...,x_i,...,x_n)$ because $f(x_1,...,x_i,...,x_n)$ is not contained in at least one part of the domain, leading to an undefined result. Thus, $f$ cannot be an Ouroboros function, which contradicts our initial assumption, so Lemma 3 holds for every $A_i$. $\qedsymbol$ 

\section*{Conclusion}

We have constructed an important type of function space that consists of all of the functions that satisfy the enigmatic Ouroboros equation \cite{1}. We observed that this notion can be extended to functions of multiple variables, and gave several examples of higher dimensional Ouroboros functions. In this exploration, we discovered that we can write the arithmetic average as a function of $n$ variables, and, more interestingly, that this arithmetic average function is indeed an Ouroboros function. Finally, we made use of the Strong Law of Large Numbers to prove that the expected value of a random variable can almost surely be an Ouroboros function for an infinite domain under certain conditions. We also pointed out that there can be Ouroboros functions for mixed domains. The concept of a self-referential function is imbued with immense mathematical, scientific, and philosophical significance, and we have linked the concept of an Ouroboros space in functional analysis to probability theory. Undoubtedly, further investigation into these concepts will yield indescribably fruitful results going forward.


\begin{thebibliography}{9}

\bibitem{1} Soto-Andrade, J., Jaramillo-Riveri, S., Gutiérrez, C., \& Letelier, J. (2011). Ouroboros avatars: A mathematical exploration of self-reference and metabolic closure. ECAL.

\bibitem{2}
Evans, L. C. (2013). 2.5. Law of Large Numbers, Central Limit Theorem. In An Introduction to Stochastic Differential Equations (pp. 24–25). American Mathematical Society. 
\end{thebibliography}
\end{document}